\documentclass[letterpaper, 10 pt, conference]{ieeeconf}
\IEEEoverridecommandlockouts
\usepackage{amsmath}
\PassOptionsToPackage{hyphens}{url}\usepackage[hidelinks]{hyperref}
\usepackage{amssymb}
\usepackage{graphicx}
\usepackage[font=small]{caption}
\usepackage{multirow}
\usepackage{makecell}
\usepackage{float}
\usepackage{xcolor}

\title{\LARGE \bf
Scaling Mixed-Integer Programming for Certification of Neural Network Controllers Using Bounds Tightening
}

\author{Philip Sosnin and Calvin Tsay
\thanks{The authors are with the Department of Computing, Imperial College London, UK. Email: {\tt\small p.sosnin23@imperial.ac.uk}, {\tt\small c.tsay@imperial.ac.uk}}
}

\begin{document}

\maketitle
\thispagestyle{empty}
\pagestyle{empty}

\begin{abstract}

Neural networks offer a computationally efficient approximation of model predictive control, but they lack guarantees on the resulting controlled system’s properties. Formal certification of neural networks is crucial for ensuring safety, particularly in safety-critical domains such as autonomous vehicles. One approach to formally certify properties of neural networks is to solve a mixed-integer program based on the network. This approach suffers from scalability issues due to the complexity of solving the resulting mixed-integer programs. 
Nevertheless, these issues can be (partially) mitigated via bound-tightening techniques prior to forming the mixed-integer program, which results in tighter formulations and faster optimization. This paper presents bound-tightening techniques in the context of neural network explicit control policies. Bound tightening is particularly important when considering problems spanning multiple time steps of a controlled system, as the bounds must be propagated through the problem depth. Several strategies for bound tightening are evaluated in terms of both computational complexity and tightness of the bounds.

\end{abstract}

\section{Introduction}

Model predictive control (MPC) is an optimal control scheme that has gained popularity due to its flexibility in handling constrained nonlinear systems in a wide variety of settings. MPC computes the optimal control action by solving an optimization problem on a forecast of the system over a limited horizon \cite{rawlings_model_2017}. This optimization problem takes the current state of the system as an initial condition and can include arbitrary system dynamics and state constraints. In its implicit form, this optimization problem is solved online in real time at each control interval. This has traditionally limited the application of MPC to either small systems, simplified models, or systems with long control intervals due to the time required to solve the MPC problem.

The MPC optimization problem for linear time-invariant systems subject to affine constraints and a quadratic cost function takes the form of a quadratic program. This allows for the solution set, which can be represented as a piecewise-affine (PWA) function, to be computed offline using multi-parametric programming and applied to the system---a paradigm known as explicit MPC. Although this moves much of the heavy computation offline, the number of regions in the PWA control law grows exponentially with the dimension of the system and the prediction horizon \cite{alessio_survey_2009}. Therefore, there can still be significant memory requirements and computational costs associated.

Current work on overcoming the limitations of explicit MPC can be roughly grouped into two categories.
The first is reducing the computational and memory requirements through more efficient representations of the PWA controller \cite{tondel_evaluation_2003}. The second is \textit{approximate} explicit MPC schemes, which can reduce the memory and computational requirements at the cost of a sub-optimal control scheme \cite{parisini_receding-horizon_1995}. Deep neural networks (NNs) are particularly promising as explicit control laws, as they allow for a large number of PWA regions to be represented with relatively few parameters. NNs come with a large offline training cost, but the inference cost and memory requirements are low, allowing them to be deployed even on embedded systems \cite{karg_efficient_2020}.

Although standard MPC schemes provide performance, feasibility, and stability guarantees, these are lost when an approximate scheme is used. In the context of NN-based approximations, recent works show that formal guarantees may be recovered by analysing properties of the NN model. For example, closed-loop stability can be guaranteed by checking if the worst-case approximation error of the NN falls within the disturbance set of a robust MPC scheme \cite{schwan_stability_2023}. Another approach is to combine the worst case error of the NN with its Lipschitz constant to provide sufficient conditions for the stability of the system \cite{fabiani_reliably-stabilizing_2023}. Alternatively, it has been shown that the linear system under NN control can be made stable if the reachable set lies within some safe subset of the state space \cite{karg_stability_2020}.

Mixed-integer programming (MIP) provides a framework in which NN controller properties such as the above can be certified. Specifically, NNs with rectified linear (ReLU) activation functions can be encoded exactly using mixed-integer linear constraints~\cite{huchette2023deep}, allowing optimization problems such as computing the worst-case approximation error to be solved to global optimality using MIP.
Unfortunately, these problems suffer from scalability issues due to the complexity of MIP, but these issues can be partially mitigated by ensuring that MIP formulations with tight linear relaxations are used.
MIP formulations of NNs require valid bounds on the intermediate variables in the network, such as the neuron activations. The tightness of these bounds directly influences the tightness of the formulation and subsequently, the speed of the MIP solver. 

In this work, we consider several strategies for obtaining tight intermediate bounds in the context of certifying NN controllers. We first investigate how bounds are propagated through both NN layers and dynamic linear systems.  We then computationally evaluate various strategies, including interval and optimization-based methods. For the latter, we present strategies based on both the standard big-M, as well as more recent, stronger formulations. 

\section{Preliminaries}
\subsection{Model predictive control (MPC)}
Model predictive control is an optimization-based control scheme that relies on a forecast of a dynamical system over a limited horizon. The system is often approximated as a discrete linear time-invariant system of the form:
\begin{equation} \label{eq:lti}
x_{k+1} = Ax_{k} + Bu_{k}
\end{equation}
where $x_{k} \in \mathbb{R}^{n_x}$ is the state vector, $u_{k} \in \mathbb{R}^{n_u}$ is the control input vector, $A \in \mathbb{R}^{n_x \times n_x}$ is the matrix of system dynamics, and $B \in \mathbb{R}^{n_x \times n_u}$ is the matrix of control dynamics.
Using a quadratic cost function, the MPC control law is obtained by solving the following at every control interval \cite{rawlings_model_2017}:
\begin{subequations}\label{eq:implicit_mpc}
\begin{align}
\underset{x, u}{\text{min }} & x_N^T P x_N+\sum\limits_{k=0}^{N-1} x_k^T Q x_k + u_k^T R u_k \\
\text { s.t. } & x_0=x_{\text {init}},  \label{eq:implicit_mpc_b}\\
& x_{k+1}=A x_k+B u_k, \quad \forall k \in [N-1] \label{eq:implicit_mpc_c}\\
& u_k \in \mathcal{U}, \hspace{2.2cm} \forall k \in [N-1] \label{eq:implicit_mpc_d}\\
& x_k \in \mathcal{X}, \hspace{2.2cm} \forall k \in [N] \label{eq:implicit_mpc_e}
\end{align}
\end{subequations}
where $P \in \mathbb{R}^{n_x \times n_x}, Q \in \mathbb{R}^{n_x \times n_x}$ and $R \in \mathbb{R}^{n_u \times n_u}$ are the terminal, state, and control cost matrices. The state and control constraint sets $\mathcal{X}$ and $\mathcal{U}$ are assumed to be bounded polytopic sets given by $\{x \in \mathbb{R}^{n_x} \mid C_x x \leq c_x\}$ and $\{u \in \mathbb{R}^{n_u} \mid C_u u \leq c_u\}$, for some $C_x \in \mathbb{R}^{n_{cx} \times n_x}$ and $C_u \in \mathbb{R}^{n_{cu} \times n_u}$. Starting from an initial condition $x_0 = x_\text{init}$, the optimal solution $u^*_0$ is applied as the control input to the system. We denote this optimal solution for an arbitrary initial condition $x$ to be  $u_\mathit{MPC}(x)$. 
Note that \eqref{eq:implicit_mpc_c} denotes the dynamic model of the type \eqref{eq:lti} and could be replaced with, e.g., a nonlinear model for increased accuracy. 

\subsection{Explicit MPC using neural networks (NNs)}

Solving problem $\eqref{eq:implicit_mpc}$ in real-time can be computationally prohibitive, particularly on embedded systems. 
This is especially true when a nonlinear/discrete model of the system is used (termed nonlinear/hybrid MPC respectively). 
An alternative is explicit MPC, where the solution set to $\eqref{eq:implicit_mpc}$ can be pre-computed offline, e.g., using multi-parametric programming \cite{alessio_survey_2009}. 
Nevertheless, the previously mentioned computational limitations of both implicit and explicit MPC have motivated approximate control schemes. One promising avenue towards this is to approximate the MPC control law using a feedforward NN model, given by:
\begin{alignat}{2}
\hat{z}^{(l)} &= W^{(l)}z^{(l-1)} + b^{(l)}, \qquad &&\forall l \in [L]\\
z^{(l)} &= \sigma\left(\hat{z}^{(l)}\right), &&\forall l \in [L-1]
\end{alignat}
where $W^{(l)} \in \mathbb{R}^{n_l \times n_{l-1}}$ and $b^{(l)} \in \mathbb{R}^{n_l}$ are the parameters of the $l$-th layer of the network and $\sigma$ is the activation function which we take to be the rectified linear unit (ReLU) for all neurons. The input to the network is $z^{(0)}$, and the output is $f_\mathit{NN}\left(z^{(0)}\right)=\hat{z}^{(L)}$.

Although more advanced NN architectures, e.g., residual or recurrent networks, are possible, simple feedforward ReLU networks are sufficient to exactly represent an MPC law \cite{karg_efficient_2020}. Furthermore, \cite{karg_efficient_2020} establishes theoretical limits on the minimum NN sizes needed to accurately approximate an MPC controller, which are typically compact ($<$100 neurons) for low-dimensional LTI systems.

To train an NN model to approximate an MPC controller, training data are obtained by choosing $n_\mathit{init}$ random initial conditions $x_\mathrm{init}$ and generating trajectories by simulating $n_\mathit{steps}$ of the closed-loop problem under the MPC control law $\eqref{eq:implicit_mpc}$. This yields $n_\mathit{train} = n_\mathit{init} \times n_\mathit{steps}$ training pairs $(x_i, u_\mathit{MPC}(x_i))$ which can be used to train the parameters of the network to minimize mean squared error using a standard optimizer. 
Note that we consider control intervals of one step in this work, but larger control intervals are easily handled by appending the dimensionality of $u_\mathit{MPC}(x_i)$.

\subsection{Mixed-integer programming (MIP) formulations}

In order to obtain formal guarantees on the behaviour of an NN controller, we wish to solve optimization problems over NNs to global optimality. In particular, ReLU NNs can be represented exactly using MIP~\cite{huchette2023deep}. Additionally, the solutions sets of parametric quadratic programs (including MPC policies) can be formulated exactly as MIPs \cite{schwan_stability_2023}. As such, optimization problems that include both NNs and MPC policies can be formulated as combined MIPs and solved to global optimality using off-the-shelf MIP solvers.

\textbf{Big-M formulation for ReLU.} A simple way to represent the ReLU activation function, $\operatorname{max}\left(W_i^{(l)}z^{(l-1)} + b_i^{(l)}, 0\right)$, is to use the so-called `Big-M' formulation \cite{lomuscio_approach_2017}:
\begin{subequations}\label{eq:bigm_formulation}
\begin{align}
&z^{(l)}_i \geq W_i^{(l)}z^{(l-1)} + b_i^{(l)} \\
&z^{(l)}_i \leq\left(W_i^{(l)}z^{(l-1)} + b_i^{(l)}\right)-\underline{M}_{i}^{(l)}(1-s) \\
&0 \leq\  z^{(l)}_i \leq \overline{M}_{i}^{(l)} s \\
& s \in\{0,1\}
\end{align}
\end{subequations}
where $W^{(l)}_i$ is the $i$-th row of the weight matrix, $s$ is an auxiliary binary variable denoting the on/off state of the neuron, and $\underline{M}_{i}^{(l)},  \overline{M}_{i}^{(l)}$ are `Big-M constants' that satisfy 
\begin{align} \label{eq:bigmbounds}
\underline{M}_{i}^{(l)} \leq \min_{{z}^{(l-1)}} W_i^{(l)}z^{(l-1)} + b_i^{(l)} \leq 
\overline{M}_{i}^{(l)}
\end{align}
From this we can observe that the big-M constants must be valid bounds of the pre-activation term. 
We call any function with such a mixed-integer formulation `MIP-representable.' Since the composition of MIP-representable functions is also MIP-representable \cite{schwan_stability_2023}, we have that $f_{NN}\left(z^{(0)}\right)=\hat{z}^{(L)}$ can be represented in mixed-integer form by adding constraints $\eqref{eq:bigm_formulation}$ for every neuron in the network.

\textbf{Partitioned formulation for ReLU.} While the big-M formulation represents an NN exactly over the box domain \eqref{eq:bigmbounds}, it suffers from a weak convex relaxation \cite{anderson_strong_2019}. As discussed in more detail in Section $\ref{sec:bounds}$, a weak convex relaxation can lead to longer solution times when using a branch-and-bound solver. Stronger formulations, such as the extended formulation \cite{anderson_strong_2019}, can represent the convex hull of each neuron exactly using $\mathcal{O}(n_{l-1})$ auxiliary variables and constraints for each neuron in the network. As such, the extended formulation results in larger problem sizes that can lead to longer solve times despite its tighter relaxation.

Our recent work introduces an intermediate formulation that can balance the trade-off between tightness and compactness \cite{tsay_partition-based_2021}. Specifically, the sum $W_i^{(l)}z^{(l-1)}$ is partitioned into several auxiliary variables $\tilde{z}_p$ for $p=1, \dots, P$ such that
\begin{subequations}\label{eq:partitioned}
\begin{align}
W_i^{(l)}z^{(l-1)} &= \sum\limits_{j=1}^{n_{l-1}}W_{ij}^{(l)}z^{(l-1)}_j = \sum\limits_{p=1}^P \tilde{z}_p\\
\tilde{z}_p &= \sum\limits_{j \in \mathbb{S}_p} W_{ij}^{(l)}z^{(l-1)}_j
\end{align}
\end{subequations}
where the partitions $\mathbb{S}_p$ are chosen such that $\mathbb{S}_1 \cup \dots \cup \mathbb{S}_P = \{1, \dots, n_{l-1}\}$ and $\mathbb{S}_j \cap \mathbb{S}_i= \emptyset$, $\forall i \neq j$.
Strategies for choosing the partitions are described in \cite{tsay_partition-based_2021}. 
An extended formulation is then applied in the space of the $\tilde{z}^{(l)}_i$ variables, resulting in a potentially tighter formulation.
Specifically, in cases where $P=1$, the partitioned formulation shares an equivalent convex relaxation with the big-M formulation.
Conversely, when $P=n_{l-1}$, the partitioned formulation aligns with the full extended formulation.
Varying $P$ between these two extremes enables a balance to be struck between the tightness and size of the intermediate formulations.

Note that, similar to the big-M formulation, intermediate bounds on each partitioned variable $\tilde{z}_p$ are required for the partitioned formulation. For both formulations \eqref{eq:bigm_formulation} and \eqref{eq:partitioned}, the tightness of the involved intermediate bounds directly influences the tightness of the overall convex relaxation.

\textbf{Mixed-integer formulation of implicit MPC.} 
The linear MPC control problem $\eqref{eq:implicit_mpc}$ can be encoded in a MIP by considering the KKT optimality conditions \cite{schwan_stability_2023}. In particular, the stationarity, primal feasbility, and dual feasibility conditions can be encoded with linear constraints, while the complementarity conditions result in bilinear constraint(s). We note that these bilinear constraints can be linearized through the use of big-M constants.

This MIP encoding is a potentially larger and more costly problem than the original quadratic program $\eqref{eq:implicit_mpc}$ but it allows for the MPC control law to be evaluated by solving a feasibility problem rather than an optimization problem. This enables the MPC controller to be embedded into larger MIPs. 

\subsection{Certification problem statement}

Recent works employ MIP to certify properties of NN control policies, such as the worst case error \cite{schwan_stability_2023}, the Lipshitz constant \cite{fabiani_reliably-stabilizing_2023}, reachable sets \cite{lomuscio_approach_2017}, and closed-loop stability \cite{karg_stability_2020}.
Although other ``complete'' NN certification techniques such as SMT~\cite{katz2017reluplexefficientsmtsolver} exist, we focus on MIP as a versatile framework for analyzing the properties of approximate MPC (e.g., the above encoding for MPC schemes).
In this work, we consider the two following certification problems.

\textbf{Worst-case error.} The worst-case approximation error of the NN controller can be calculated using the following MIP:
\begin{subequations}\label{eq:worst_case}
\begin{align}
\max\limits_{x_0} \ & \|u - \hat{u} \| \label{eq:worst_case_a}\\
\text{ s.t. } & x_0 \in \mathcal{X}_{in} \label{eq:worst_case_b}\\
& u = u_\mathit{MPC}(x_0) \label{eq:worst_case_c}\\
& \hat{u} = f_\mathit{NN}(x_0) \label{eq:worst_case_d}
\end{align}
\end{subequations}
where \eqref{eq:worst_case_a} and \eqref{eq:worst_case_b} contain, respectively, a MIP-representable norm and input domain. 
Constraint \eqref{eq:worst_case_c} encodes the MPC model, e.g., using the KKT conditions, and \eqref{eq:worst_case_d} encodes the NN controller, e.g., using the Big-M, partitioned or extended formulations. 
The result of \eqref{eq:worst_case} can additionally be used to certify the stability of approximate control schemes \cite{schwan_stability_2023}, \cite{fabiani_reliably-stabilizing_2023}.

\textbf{Reachability analysis over multiple time-steps.} Another useful quantity is the reachable set, which can be bounded using MIP \cite{karg_stability_2020}. We define the $K$-step reachable set to be the set of possible final states $x_K$ under the NN control policy $x_k = Ax_{k-1} + Bf_{NN}(x_{k-1})$ given an initial state set $x_0 \in \mathcal{X}_{in}$. We choose to bound this set with a polytope $\{x_K \mid C_{out} x_K \leq c_{out}\}$, where the matrix $C_{out}$ is chosen ahead of time and defines the facets of our polyhedral bounding set. We wish to obtain the tightest possible constants $c_{out}$ such that $C_{out} x_K \leq c_{out}$ holds for all $x_0 \in \mathcal{X}_{in}$. We can achieve this by solving the following mixed-integer program for each row of $C_{out}$:
\begin{subequations}\label{eq:reachability}
\begin{align}
\max\ & C_{out, i} x_K \label{eq:reachability_a}\\
\text{s.t.}\  & x_0 \in \mathcal{X}_{in} \label{eq:reachability_b} \\
& x_{k+1} = A x_k + B u_k, & k \in [K - 1] \label{eq:reachability_c}\\
& u_k = f_\mathit{NN}(x_k),  & k \in [K - 1] \label{eq:reachability_d}
\end{align}
\end{subequations}
where \eqref{eq:reachability_b} again includes a MIP-representable input domain, and  \eqref{eq:reachability_d} again includes the NN controller using a suitable MIP formulation. Each $C_{out, i} x_K \leq c_{out, i}$ defines a hyperplane bounding the reachable set. By choosing the directions and number of hyperplanes it is possible to get a more refined view of the output reachable set.

\section{Deriving and tightening bounds}\label{sec:bounds}

Both the big-M \eqref{eq:bigm_formulation} and partitioned \eqref{eq:partitioned} formulations of NNs described in the previous section require valid bounds for the intermediate variables of the network.
In general, the tightness of variable bounds can strongly influence the tightness of the convex relaxation for a MIP. 
The tightness of a formulation is characterized by the gap between the optimal mixed-integer solution and the solution to its linear/convex relaxation.
Branch and bound, the most common MIP algorithm, uses the solution to the relaxation to bound the objective value and prune the search space. 
As such, tighter formulations can significantly reduce the solve time (although larger formulations may require more time to solve the relaxed sub-problems). 

This section presents several methods for deriving valid bounds on the pre-activation variables of the network, with a focus on problems \eqref{eq:worst_case}--\eqref{eq:reachability}.
This collection of methods is non-comprehensive and a number of other NN verifiers may be suitable for the bound-tightening problem, such as domain propagation \cite{singh_abstract_2019}.

\begin{figure*}[ht]
    \centering
    \includegraphics[width=\textwidth,]{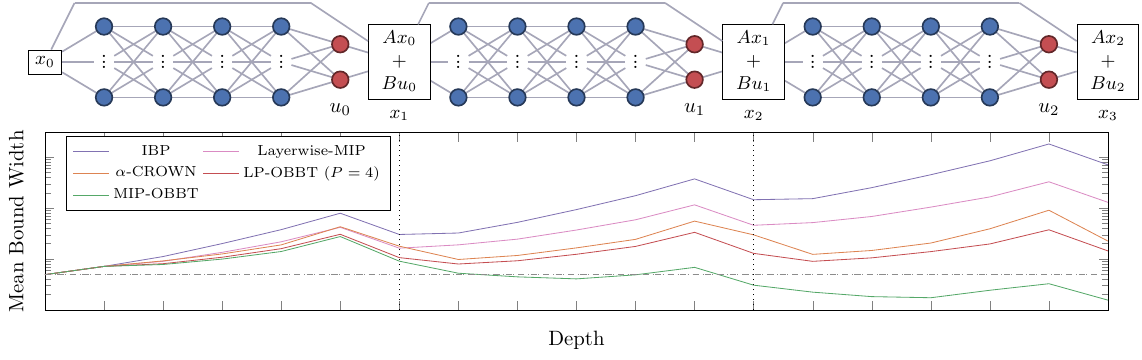}
    \caption{Comparison of bound tightening methods on a 4-layer NN controller over $K=3$ time-steps (note log-scale y-axis).}
    \label{fig:bounds_comparison}
\end{figure*}

\subsection{Interval bound propagation (IBP)}
IBP propagates an axis-aligned bound in layer-by-layer steps through the NN to obtain valid bounds for each activation variable \cite{gowal_effectiveness_2019}. In particular, given a valid bound for the $(l-1)$-th layer, ${z}^{(l-1)} \in \left[l^{(l-1)}, u^{(l-1)}\right]$, IBP computes valid bounds for $z^{(l)} = \operatorname{ReLU}\left(W^{(l)}z^{(l-1)} + b^{(l)}\right)$:
\begin{align*}
\mu = \frac{u^{(l-1)} + l^{(l-1)}}{2}, \qquad r = \frac{u^{(l-1)} - l^{(l-1)}}{2} \\
l^{(l)} = \operatorname{ReLU}\left(W^{(l)}\mu - \left|W^{(l)}\right|r + b^{(l)} \right) \leq z^{(l)}\\
u^{(l)} = \operatorname{ReLU}\left( W^{(l)}\mu + \left|W^{(l)}\right|r + b^{(l)} \right) \geq z^{(l)}
\end{align*}
IBP is used in the majority of MIP studies owing to its simplicity, but can produce loose bounds as dependencies between neurons are omitted. In fact, the bounds produced by IBP exponentially weaken as they are propagated through the depth of the NN \cite{badilla_computational_2023}.

\subsection{Linear bound propagation (CROWN)}
Instead of propagating an interval for each intermediate variable, it is possible to propagate symbolic linear bounds. This allows for dependencies between variables in the network to be captured and tighter bounds to be obtained. In order for the propagation procedure to be computationally efficient, typically only a single upper and lower bound are used to bound each ReLU activation function \cite{singh_abstract_2019}. 
Specifically, the upper bound matches that of the big-M formulation \eqref{eq:bigm_formulation}, while the lower bound interpolates between the constraints $z^{(l)}\geq 0$ and $z^{(l)} \geq \hat{z}^{(l)}$. 

CROWN \cite{zhang_efficient_2018} adaptively chooses between these two linear lower bounds, while other methods employ a relaxation with parallel upper and lower bounds. It has been shown that all such linear bound propagation algorithms are equivalent to projecting out intermediate network layers using Fourier-Motzkin elimination \cite{tjandraatmadja_convex_2020}.

The later $\alpha$-CROWN \cite{xu_fast_2021} employs gradient descent to jointly optimize the choice of the lower bound slopes for all neurons in the network. 
Specifically, the lower bound is defined using a slope $\alpha \in [0,1]$ such that $z^{(l)} \geq \alpha \hat{z}^{(l)}$, and the values of $\alpha$ for all neurons are optimized. 
When fully converged, this method produces bounds that are equivalent to the relaxed Big-M based OBBT described below.

\subsection{Optimization-based bounds tightening (OBBT)}
As implied by its name, OBBT solves optimization problems to compute the minimum and maximum values reachable by each variable. For instance, if the original MIP formulation is used, OBBT can compute the exact tightest bounds, but at the cost of solving 2 MIPs for every intermediate variable. Note that when obtaining bounds for an intermediate variable, only those variables in prior layers of the problem must be included in the OBBT subproblem. This can still be prohibitively expensive, particularly for large problems defined over multiple timesteps, where later OBBT subproblems may have a complexity similar to the full certification problem. 
We refer to this (exact) bound tightening technique as MIP-OBBT.
The Big-M and partitioned formulations of ReLU neural networks are both exact, so MIP-OBBT yields the same solution regardless of which formulation is used.

The computational cost of MIP-OBBT can be greatly reduced by instead considering a linear relaxation of the OBBT subproblems, which typically results in looser bounds.
The linear relaxation is formed by relaxing binary variables to be continuous in $[0, 1]$. 
We will refer to this as linear programming-based OBBT (LP-OBBT).
Unlike MIP-OBBT, the choice of formulation directly influences the tightness of the resulting bounds.
As discussed above, the Big-M formulation has a particularly weak linear relaxation so the intermediate bounds obtained from LP-OBBT may be a large over-approximation of the true bounds.
On the other hand, partitioned formulations lead to tighter bounds as the number of partitions $P$ increases, at the cost of more expensive subproblems. 

\subsection{Layer-wise OBBT}

Instead of solving the entire formulation for each intermediate variable, we can decompose the problem into smaller subproblems that will still give valid, but weaker bounds. For example, we can form layer-wise subproblems that compute the min/max of the output of a layer in the NN given bounds on the input to the layer. Solving this smaller (MI)LP subproblem accounts for dependencies in single layers of the network and therefore we can obtain relatively tight bounds while keeping the subproblems small. By starting with the first layer and forming subproblems in a rolling-horizon fashion, we can obtain valid bounds for the entire problem \cite{zhao_bound_2024}. For deep networks, the layer-wise subproblems are significantly smaller than the overall MIP formulation.
For a network of fixed width, the time-complexity of layer-wise OBBT scales linearly with the depth, while full OBBT scales exponentially.
In general, more than one layer can be considered per rolling-horizon subproblem, providing a trade-off between speed and tightness of the bounds.
In this paper, we only consider the rolling horizon procedure with a single layer per subproblem.

\subsection{Combining bound tightening methods}
The various bound-tightening techniques outlined above can be combined to achieve more flexible trade-offs between solution tightness and computational time.
For instance, starting with linear bound propagation can generate tighter initial formulations that accelerate subsequent OBBT subproblem solution times.
Additionally, different strategies can be applied selectively across the network; for example, using tighter methods/formulations in the early layers and faster, less precise methods in the deeper layers. 
Due to space limitations, this work focuses on analyzing each method in isolation to evaluate its individual trade-off characteristics. 
Identifying the optimal combination of these techniques for specific neural network architectures remains an open area for future work.

\section{Bounding multi-timestep problems}

The bound-tightening methods described above are typically applied to feedforward NNs. In the case of problem \eqref{eq:reachability}, we wish to bound variables over multiple time steps of the linear system, including multiple instances of the same NN controller. As the system dynamics in \eqref{eq:lti} are linear, the state update can be viewed as an additional affine layer in an augmented ``network.'' This allows the entire multi-step problem to be unrolled into a single deep NN, with skip connections carrying each state forward into the next timestep of the problem. This emphasizes the need for tight bounding methods, as even shallow NNs result in extremely deep problems when multiple time steps are considered.

Figure \ref{fig:bounds_comparison} shows the evolution of different bounding methods throughout a multi-timestep problem. If the network is well trained, it is expected that the state of the dynamical system should be brought back to its set point. As such, the `true' bounds on the states of the linear system should contract with each time step. MIP-OBBT produces exact bounds on each variable and as such is the only method which exhibits shrinking bounds with each time step. The over-approximation introduced by each other method results in bounds that grow exponentially with the depth of the problem, but it should be noted that the tighter methods all exhibit significantly slower growth than IBP.

\begin{table*}[h]
\begin{center}
\begin{tabular}{|cl|cccccccc|}
\cline{3-10}
\multicolumn{2}{c|}{} & IBP & CROWN & $\alpha$-CROWN  & \makecell{LP-OBBT\\(Big-M)} & \makecell{LP-OBBT\\(P=2)} & \makecell{LP-OBBT\\(P=4)} & \makecell{Layer-wise\\OBBT} & MIP-OBBT \rule{0pt}{4ex} \\
\hline
\multirow{4}{*}{6x10} & Avg. Bound Width & 130.29 & 51.64 & 36.53 & 33.37 & 27.61 & 26.52 & 33.11 & \textbf{22.78} \rule{0pt}{2ex} \\
 & Avg. Bound Time (s) & \textbf{0.04} & 0.02 & 2.89 & 0.67 & 2.10 & 3.35 & 1.22 & 5.47 \\
 & Avg. Solve Time (s) & 1238.26 & 1098.76 & 1189.00 & 1081.02 & 1016.70 & 982.35 & 935.87 & \textbf{884.71}  \\
 & Number Solved & 67 & 66 & 63 & 67 & 69 & 68 & 68 & \textbf{73} \\
 \hline
\multirow{4}{*}{1x60} & Avg. Bound Width & 32.56 & 30.83 & 30.35 & 30.35 & 29.02 & \textbf{28.76} & - & \textbf{28.76} \rule{0pt}{2ex} \\
 & Avg. Bound Time (s) & \textbf{0.0004} & 0.05 & 0.53 & 0.31 & 1.12 & 2.04 & - & 0.26 \\
 & Avg. Solve Time (s) & 2172.99 & 2051.31 & 2729.05 & \textbf{1795.27} & 1802.75 & 1815.76 & - & 1820.14 \\
 & Number Solved & 19 & 19 & 25 & 23 & 24 & 24 & - & \textbf{25} \\
 \hline
\end{tabular}
\end{center}
\caption{Worst case error results with an MPC horizon of $N =9$. Avg. solve time is computed only for problems solved by all methods.}
\label{tab:worst_case}
\end{table*}

\begin{table*}[h]
\begin{center}
\begin{tabular}{|cl|cccccccc|}
\cline{3-10}
\multicolumn{2}{c|}{} & IBP & CROWN & $\alpha$-CROWN  & \makecell{LP-OBBT\\(Big-M)} & \makecell{LP-OBBT\\(P=2)} & \makecell{LP-OBBT\\(P=4)} & \makecell{Layer-wise\\OBBT} & MIP-OBBT \rule{0pt}{4ex} \\
\hline
\multirow{4}{*}{6x10} & Avg. Relaxation Gap (\%) & 2217423.69 & 463991.68 & 8302.14 & 5015.85 & 4168.03 & 3134.31 & 66959.54 & \textbf{40.20} \rule{0pt}{2ex} \\
 & Avg. Bound Time (s) & \textbf{0.03} & 1.35 & 65.20 & 9.36 & 38.08 & 85.23 & 1.65 & 1125.94 \\
 & Avg. Solve Time (s) & 996.27 & 933.41 & 571.76 & 423.15 & 344.89 & 271.53 & 809.08 & \textbf{5.72} \\
 & Number Solved & 76 & 73 & 82 & 86 & 86 & 96 & 76 & \textbf{100} \\
 \hline
 \multirow{4}{*}{1x60} & Avg. Relaxation Gap (\%) & 15088.99 & 849.54 & 337.69 & 292.72 & 236.22 & 195.31 & 12926.87 & \textbf{24.54}  \rule{0pt}{2ex} \\
 & Avg. Bound Time (s) & \textbf{0.03} & 0.53 & 4.01 & 3.57 & 16.16 & 36.39 & 0.52 & 431.36 \\
 & Avg. Solve Time (s) & 62.57 & 75.81 & 49.68 & 43.62 & 32.64 & 29.60 & 65.01 & \textbf{7.80} \\
 & Number Solved & 100 & 100 & 100 & 100 & 100 & 100 & 100 & 100 \\
 \hline
\end{tabular}
\end{center}
\caption{Reachability results for $K=4$ timesteps. Avg. solve time is computed only for problems solved by all methods.}
\label{tab:reachability}
\end{table*}

\section{Case study}

We consider the 4-dimensional linear system from \cite{jones2010polytopic} with the following state and control dynamics
\begin{align*}
x_{k+1}=\left[\begin{array}{cccc}
0.7 & -0.1 & 0.0 & 0.0 \\
0.2 & -0.5 & 0.1 & 0.0 \\
0.0 & 0.1 & 0.1 & 0.0 \\
0.5 & 0.0 & 0.5 & 0.5
\end{array}\right] x_k+\left[\begin{array}{cc}
0.0 & 0.1 \\
0.1 & 1.0 \\
0.1 & 0.0 \\
0.0 & 0.0
\end{array}\right] u_k
\end{align*}
with constraints from \cite{chen_approximating_2018} given by
\begin{align*}
\left|x_k\right| \leq \left[6.0, 6.0, 1.0, 0.5\right]^\top, \quad \left|u_k\right| \leq\left[5.0, 5.0\right]^\top.
\end{align*}
The state and control cost matrices are given by $Q=I^{4\times 4}$ and $R=I^{2 \times 2}$. The terminal cost matrix $P$ is taken to be the solution of the discrete-time Riccati equation \cite{fink_implementation_2021}. 
To ensure recursive feasibility, the input domain for the MPC problem is restricted to the maximum control-invariant set computed using MPT3 \cite{herceg_multi-parametric_2013}. 
This maximum invariant set is further used as $\mathcal{X}_{\text{in}}$ for both the worst case error and reachability analysis problems.
Two simple feedforward NN controllers are considered: a `deep' six-layer network with 10 neurons per layer and a `wide' network with one layer of 60 neurons.
Each network was trained using $n_{\text{train}} = 30000$ samples generated from the true MPC controller with a horizon of $N = 9$. 
Intermediate bounds are first computed using a variety of methods: IBP, ($\alpha$-)CROWN, LP-OBBT, MIP-OBBT, and layer-wise OBBT. These bounds are then used to formulate certification problems, and we investigate the resulting solve times for each formulation (within a lime limit of 3600s).
All experiments are run on a server equipped with AMD EPYC 7742 CPUs (2.25 GHz) with a maximum of 32 threads and solved using Gurobi 10 \cite{gurobi_optimization_llc_gurobi_2023}.

The computations for IBP and ($\alpha$-)CROWN are carried out using the auto-LiRPA library \cite{xu2020automatic}.
While ($\alpha$-)CROWN typically benefits from GPU acceleration, the performance difference for NNs of this size is minimal and, in some cases, can be hindered by memory overheads. $\alpha$-CROWN is executed for a fixed 20 iterations.
LP-OBBT, MIP-OBBT, and Layer-wise OBBT are computed using a custom implementation (\url{www.github.com/psosnin/Bounds-Tightening-for-Neural-Network-Control}). 

\subsection{Worst-case error results}

In this section, we consider problem \eqref{eq:worst_case}, which calculates the worst-case error between an NN approximation and the true MPC controller. Here we encode the NN using the Big-M constraints \eqref{eq:bigm_formulation} and constrain the initial condition to be within the maximum control invariant set such that the MPC problem is feasible everywhere in the input domain. 100 NNs were trained on distinct MPC datasets and intermediate bounds calculated using each of the bound tightening methods described above.

The results for the worst case error problem instances are shown in Table \ref{tab:worst_case}. 
As the linear relaxation of problem \eqref{eq:worst_case} is unbounded, we report the tightness of the bounds by the average width of the bounds of the neural network output.
For the deep network, formulations with tighter bounds generally result in faster solve times.
The effect is less significant for the wide network, as even the weak bounding methods are fairly tight for a single layer network.
In both cases however, the bound tightening time is significantly less than the resulting solve time.
This suggests that much of the complexity is dominated by the MIP formulation of the MPC controller, which is not affected by tightening of the NN bounds.
In turn, this suggests that for the worst case error problem it is always desirable to use the tightest possible bound tightening method.

\subsection{Reachability analysis results}

We now consider bounding the reachable set of the system under NN control \eqref{eq:reachability} over $K=4$ timesteps.
As above, 100 NN controllers of each architecture were generated, and intermediate bounds for the whole multi-step problem were computed for each network.
A hyperplane with which to bound the reachable set was randomly generated to form the objective of problem \eqref{eq:reachability}.
Although here we are only bounding the reachable set with a single hyperplane, in general this problem must be solved several times with multiple different objectives to form a bounded polyhedral over-approximation of the reachable set.
However, the intermediate bounds only need to be computed once and can be reused for each subsequent tightening direction.

Table \ref{tab:reachability} shows the results for the reachability problems. Here we see a stronger impact of bounds tightening, with methods that produce a tighter relaxation having significantly faster solve times.
This is due to the increased depth of the problem over multiple time steps, combined with the fact that the complexity of the MIP \eqref{eq:reachability} is dominated by the involved NNs.
Tighter bounds allow not only for more efficient pruning of the state space in the branch-and-bound algorithm but also for the potential `stabilisation' of neuron activations to be either on or off.

We note that the MIP-OBBT bounds form an exact axis-aligned bound of the reachable set, which may already be sufficient to certify the safety of the controller. 
Solving \eqref{eq:reachability} in this case simply amounts to computing some new bounding facet of the reachable set. 
As such, the solve times are extremely fast but at the cost of a long bound-tightening times.
Since bound tightening only needs to be performed once, but the reachability problem needs to be solved multiple times, longer bounding times might be acceptable in cases where more refined bounds of the output reachable set are required. 
Note that the OBBT problems are also easily parallelisable. 
For both architectures, the LP-OBBT methods provided a good trade-off between bound-tightening and solve times.
Although $\alpha$-CROWN can produce bounds equivalent to big-M LP-OBBT, for deep problems it requires more time to converge than simply solving the LP-OBBT problems.

\section{Conclusions}

MIP can be used to recover some formal MPC guarantees when NN approximations are used. 
This work highlighted the importance of bound tightening strategies in providing these certifications.
We presented several different bound tightening methods and their applications to multi-timestep problems.
An extensive computational study revealed that tight bounding methods can enable MIP-based NN certification techniques to scale to larger, multi-timestep problems.

\bibliographystyle{ieeeconf}
\bibliography{references}

\end{document}